\documentclass[12pt]{article}
\usepackage{graphicx}
\usepackage{graphics,amsmath, amsthm,amssymb,amsfonts}
\usepackage{pstcol,pst-node,graphics}
\newpsobject{showgrid}{psgrid}{subgriddiv=1,griddots=10,gridlabels=6pt}
\newcommand{\st}[2]{\genfrac{[}{]}{0pt}{}{#1}{#2}}
\newtheorem{theorem}{Theorem}[section]
\newtheorem{proposition}{Proposition}[section]
\newtheorem{remark}{Remark}[section]

\title{The dinner table problem: \\the rectangular case}

\author{Roberto Tauraso 
\thanks{I would like to warmly thank Alessandro Nicolosi and Giorgio Minenkov for drawing my attention to this problem.}\\
\small Dipartimento di Matematica\\[-0.8ex]
\small Universit\`a di Roma ``Tor Vergata''\\[-0.8ex]
\small 00133 Roma, Italy\\[-0.8ex]
\small \texttt{tauraso@mat.uniroma2.it}}

\begin{document}
\maketitle

\begin{section}{Introduction}

\noindent Assume that $8$ people are seated around a table 
and we want to enumerate the number of ways that they can be permuted
such that neighbors are no more neighbors after the rearrangement.
Of course the answer depends on the {\sl topology} of the table:
if the table is a circle then it is easy to check by a simple computer program 
that the permutations that verify this property are $2832$, if it is a long bar and all the persons 
seat along one side then they are $5242$. Furthermore, if it is a rectangular table with two sides then
the rearrangements are $9512$.
The first two cases are respectively described by sequences {\sl A089222} and {\sl A002464}
of the On-Line Encyclopedia of Integer Sequences \cite{Sloane}, on the other hand the rectangular case
does not appear in the literature. 

\noindent Here is a valid rearrangement for $n=8$:
\begin{center}
\begin{pspicture}(0,-0.5)(11,2.6)
\def\tableotto{\put(0,0.05){\pspolygon(0,0)(4,0)(4,0.9)(0,0.9)}}
\def\seat{\psellipse(0,0)(0.34,0.34)}
\put(0.8,0.8){\tableotto
\put(0.5,1.45){\seat \rput(0,0){1}}
\put(1.5,1.45){\seat \rput(0,0){3}}
\put(2.5,1.45){\seat \rput(0,0){5}}
\put(3.5,1.45){\seat \rput(0,0){7}}
\put(0.5,-.45){\seat \rput(0,0){2}}
\put(1.5,-.45){\seat \rput(0,0){4}}
\put(2.5,-.45){\seat \rput(0,0){6}}
\put(3.5,-.45){\seat \rput(0,0){8}}
\uput[d](2,-0.9){First dinner}}
\put(6.2,0.8){\tableotto
\put(0.5,1.45){\seat \rput(0,0){1}}
\put(1.5,1.45){\seat \rput(0,0){6}}
\put(2.5,1.45){\seat \rput(0,0){3}}
\put(3.5,1.45){\seat \rput(0,0){7}}
\put(0.5,-.45){\seat \rput(0,0){2}}
\put(1.5,-.45){\seat \rput(0,0){8}}
\put(2.5,-.45){\seat \rput(0,0){4}}
\put(3.5,-.45){\seat \rput(0,0){5}}
\uput[d](2,-0.9){Second dinner}}
\end{pspicture}
\end{center}
whose associated permutation is
$$\pi=\left(\begin{array}{cccccccc}
1&2&3&4&5&6&7&8\\
1&2&6&8&3&4&7&5
\end{array}\right).$$
For a generic number of persons $n$ the required property can be established more formally in this way:
$$|\pi(i+2)-\pi(i)|\not= 2\quad\mbox{for $1\leq i\leq n-2$.}$$
It is interesting to note that this rearrangement problem around a table has also another 
remarkable interpretation. Consider $n$ kings to be placed on a $n\times n$ board, one in each row and column, 
in such a way that they are non-attacking with respect to these different {\sl topologies} of the board:
if we enumerate the ways on a toroidal board we find the sequence
{\sl A089222}, for a regular board we have {\sl A002464}, and finally if we divide the board
in the main four quadrants we are considering the new sequence (now recorded as {\sl A110128}). 
Here is the $8$ kings displacement that corresponds to the permutation $\pi$ introduced before:
\begin{center}
\begin{pspicture}(0,0)(6,6)
\put(0,0){\includegraphics{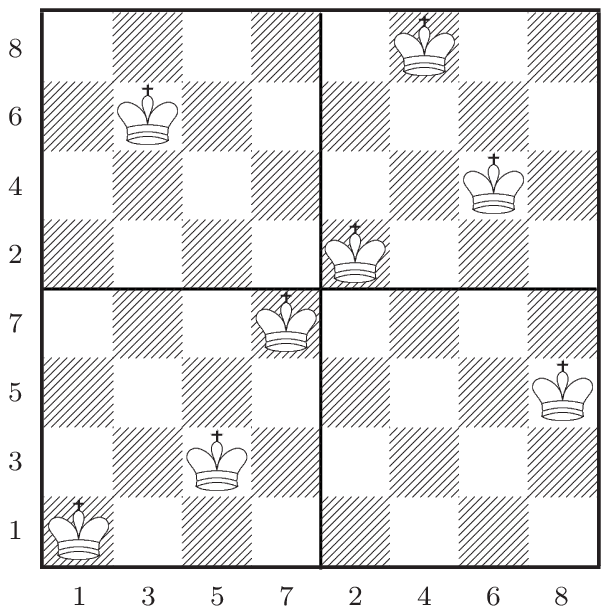}}
\end{pspicture}
\end{center}
Closed formulas for the first two sequences are known: 
for {\sl A089222} it is (see \cite{AspvallLiang})
$$a_{n,0}=\sum_{r=0}^{n-1}\sum_{c=0}^{r}(-1)^{r} 2^{c} (n-r)!\binom{r-1}{c-1}\binom{n-r}{c}\left(\frac{n}{n-r}\right)^2+ (-1)^n 2n$$
and for {\sl A002464}  it is (see \cite{Robbins})
$$a_{n,1}=\sum_{r=0}^{n-1}\sum_{c=0}^{r}(-1)^{r} 2^{c} (n-r)!\binom{r-1}{c-1}\binom{n-r}{c}$$
(remember that ${k\choose -1}=0$ if $k\not=-1$ and ${-1\choose -1}=1$).

\noindent In this paper we propose to study the sequence $a_{n,d}$ defined for $2\leq d \leq n-1$ 
as follows: $a_{n,d}$ denotes the total number of permutations $\pi$ of $\{1,2,\dots,n\}$ such that
$$|\pi(i+d)-\pi(i)|\not= d\quad\mbox{for $1\leq i\leq n-d$.}$$
The table below provides some numerical values:
\begin{center}
\begin{tabular}{| r | r | r | r | r |}
\hline
\multicolumn{1}{|c }{$n$} & 
\multicolumn{1}{|c }{$a_{n, 0}$} &
\multicolumn{1}{|c }{$a_{n, 1}$} &
\multicolumn{1}{|c }{$a_{n, 2}$} &
\multicolumn{1}{|c|}{$a_{n, 3}$} \\
\hline
\vspace{-3mm} &  &  &  &\\
1  & 1 & 1 & 1 & 1 \\
2  & 0 & 0 & 2 & 2 \\
3  & 0 & 0 & 4 & 6 \\
4  & 0 & 2 & 16 & 20\\
5  & 10 & 14 & 44 & 80\\
6  & 36 & 90 & 200 & 384 \\
7  & 322& 646 & 1288 & 2240 \\
8  & 2832& 5242 & 9512 & 15424 \\
9  & 27954& 47622 & 78652 & 123456 \\
10 & 299260& 479306 & 744360 & 1110928 \\
11 & 3474482& 5296790 & 7867148 & 11287232 \\
12 & 43546872& 63779034 & 91310696 & 127016304 \\
13 & 586722162& 831283558 & 1154292796 & 1565107248 \\
14 & 8463487844& 11661506218 & 15784573160 &  20935873872 \\
15 & 130214368530& 175203184374 & 232050062524 & 301974271248 \\
16 & 2129319003680& 2806878055610 & 3648471927912 & 4669727780624\\
\hline
\end{tabular}
\end{center}
We will show that the following formula holds for $d\geq 2$:
$$a_{n,d}=\sum_{r_1=0}^{n_1-1}\sum_{c_1=0}^{r_0}\cdots
\sum_{r_{d}=0}^{n_{d}-1}\sum_{c_{d}=0}^{r_{d}}
(-1)^{r} 2^{c} (n-r-c)!
\prod_{k=1}^{d}{n_k-r_k\choose c_k}\!\!\!\!
\sum_{\begin{array}{c}
\scriptstyle \sum_{i=1}^{c_1}l_{i,1}=r_1\\
\scriptstyle l_{i,1}\geq 1 \\
\end{array}}
\!\!\!\!\!\!\cdots\!\!\!\!\!\!
\sum_{\begin{array}{c}
\scriptstyle \sum_{i=1}^{c_{d}} l_{i,d}=r_{d}\\
\scriptstyle l_{i,d}\geq 1 
\end{array}}
\!\!\!\!\!\!\!\!q_{n,d}(L)$$
where $n_k=|N_k|=\left|\{1\leq i \leq n: i\equiv k \mod d\}\right|$, $r=\sum_{k=1}^{d} r_k$, $c=\sum_{k=1}^{d} c_k$ and
$$q_{n,d}(L)=q_{n,d}([l_1,\cdots,l_c])=
\sum_{\begin{array}{c}
\scriptstyle J_1\stackrel{\circ}{\cup}\cdots\stackrel{\circ}{\cup}J_{d}=\{1,\dots,c\}\\
\scriptstyle n_k\geq\sum_{i\in J_k}l_i \\
\end{array}}
\prod_{k=1}^{d}{n_k-\sum_{i\in J_k}l_i\choose |J_k|}\cdot|J_k|!\,.$$
In the last section, we will prove that the probability that a
permutation belongs to the set enumerated by $a_{n,d}$ tends always to $e^{-2}$ as $n$ goes to infinity. 
A more precise expansion reveals how the limit probability depends on $d$:
$${a_{n,d}\over n!}=e^{-2}\left(1+\frac{4(d-1)}{n}+O\left(\frac{1}{n^2}\right)\right).$$
\end{section}

\begin{section}{Asymptotic analysis: cases $d=0$ and $d=1$}

\begin{proposition}
The following asymptotic expansions hold:
$${a_{n,0}\over n!}\sim e^{-2}\left(1-{4\over n}+{20\over 3n^3}+{58\over 3n^4}+{736\over 15 n^5}+O\left({1\over n^6}\right)\right),$$
and
$${a_{n,1}\over n!}\sim e^{-2}\left(1-{2\over n^2}-{10\over 3n^3}-{6\over n^4}-{154\over 15 n^5}+O\left({1\over n^6}\right)\right),$$
\end{proposition}

\begin{proof} The expansion of $a_{n,0}/n!$ is contained in \cite{AspvallLiang} and it was obtained
from a recurrence relation by the method of undetermined coefficients.
However, since we are particularly interested in, we will explicitely derive the coefficient of $1/n$.

\noindent As regards $a_{n,1}/n!$, we give the detailed proof only for the coefficients of $1/n$ and $1/n^2$
(the others can be computed in a similar way).
Since 
$${a_{n,1}\over n!}=\sum_{r=0}^{n-1}\sum_{c=0}^{r}
(-1)^{r} {2^{c}\over c!}{r-1\choose c-1}{n^{\underline{r+c}}\over (n^{\underline{r}})^2}$$
and as $n$ goes to infinity the resulting series is uniformly convergent (with respect to $n$),
then we can consider each term one by one.
Moreover, it suffices to analyse the cases when $c$ is equal to $r$, $r-1$ and $r-2$ because 
the rational function $n^{\underline{r+c}}/(n^{\underline{r}})^2\sim 1/n^{r-c}$.
For $c=r$
$$
{r-1\choose r-1}{n^{\underline{2r}}\over (n^{\underline{r}})^2}\sim
{1-\st{2r}{2r-1}{1\over n}+\st{2r}{2r-2}{1\over n^2}
\over \left(1-\st{r}{r-1}{1\over n}+\st{r}{r-2}{1\over n^2}\right)^2}\sim
1-{r^{\underline{2}}+r\over n}+{r^{\underline{4}}+4r^{\underline{3}}+2r^{\underline{2}}\over 2n^2}.
$$
In a similar way, for $c=r-1$
$$
{r-1\choose r-2}{n^{\underline{2r-1}}\over (n^{\underline{r}})^2}\sim
{r-1\over n}-{(r-1)^{\underline{3}}+3(r-1)^{\underline{2}}+(r-1)\over n^2}
$$
and for $c=r-2$
$$
{r-1\choose r-3}{n^{\underline{2r-2}}\over (n^{\underline{r}})^2}\sim
{(r-2)^{\underline{2}}+2(r-2)\over 2n^2}.
$$
Hence 
\begin{eqnarray*}
{a_{n,1}\over n!}&\!\!\!\!\sim\!\!\!\!&
\sum_{r=0}^{n-1} {(-2)^{r}\over r!}\left(
1-{r^{\underline{2}}+r\over n}+{r^{\underline{4}}+4r^{\underline{3}}+2r^{\underline{2}}\over 2n^2}\right)+\\
&\!\!\!\!\!\!\!\!&-\sum_{r=1}^{n-1} {(-2)^{r-1}\over (r-1)!}\left(
{r-1\over n}-{(r-1)^{\underline{3}}+3(r-1)^{\underline{2}}+(r-1)\over n^2}\right)+\\
&\!\!\!\!\!\!\!\!&+\sum_{r=2}^{n-1} {(-2)^{r-2}\over (r-2)!}\left(
{(r-2)^{\underline{2}}+2(r-2)\over 2n^2}\right).
\end{eqnarray*}
Taking the sums we find that
\begin{eqnarray*}
{a_{n,1}\over n!}&\!\!\!\!\sim\!\!\!\!&
e^{-2}\left(1+ {(-2^2+2)+2\over n}+{(2^4-4\cdot2^3+2\cdot 2^2)+2\cdot(-2^3+3\cdot 2^2-2)+(2^2-2\cdot 2)\over 2n^2}\right)\\
&\!\!\!\!\sim\!\!\!\!& 
e^{-2}\left(1- {2\over n^2}\right).
\end{eqnarray*}
Note that the same strategy can be applied also to $a_{n,0}/n!$. For example the coefficient of $1/n$ can be easily 
calculated: for $c=r$ the formula gives 
$$
{r-1\choose r-1}{n^{\underline{2r}}\over (n^{\underline{r}})^2}\cdot\left(\frac{n}{n-r}\right)^2\sim
\left(1-{r^{\underline{2}}+r\over n}\right)\cdot\left(1+{2r\over n}\right)\sim 1-{r^{\underline{2}}-r\over n}
$$
and therefore
\begin{eqnarray*}
{a_{n,0}\over n!}&\!\!\!\!\sim\!\!\!\!&
\sum_{r=0}^{n-1} {(-2)^{r}\over r!}\left(
1-{r^{\underline{2}}-r\over n}\right)-\sum_{r=1}^{n-1} {(-2)^{r-1}\over (r-1)!}\left({r-1\over n}\right)\\
&\!\!\!\!\sim\!\!\!\!& 
e^{-2}\left(1+ {(-2^2-2)+2\over n}\right)\sim e^{-2}\left(1- {4\over n}\right).
\end{eqnarray*}
\end{proof}
\end{section}

\begin{section}{The formula for $d\geq 2$}
\begin{theorem} For $d\geq 2$
$$a_{n,d}=\sum_{r_1=0}^{n_1-1}\sum_{c_1=0}^{r_0}\cdots
\sum_{r_{d}=0}^{n_{d}-1}\sum_{c_{d}=0}^{r_{d}}
(-1)^{r} 2^{c} (n-r-c)!
\prod_{k=1}^{d}{n_k-r_k\choose c_k}\!\!\!\!
\sum_{\begin{array}{c}
\scriptstyle \sum_{i=1}^{c_1}l_{i,1}=r_1\\
\scriptstyle l_{i,1}\geq 1 \\
\end{array}}
\!\!\!\!\!\!\cdots\!\!\!\!\!\!
\sum_{\begin{array}{c}
\scriptstyle \sum_{i=1}^{c_{d}} l_{i,d}=r_{d}\\
\scriptstyle l_{i,d}\geq 1 
\end{array}}
\!\!\!\!\!\!\!\!q_{n,d}(L)$$
where $n_k=|N_k|=\left|\{1\leq i \leq n: i\equiv k \mod d\}\right|$, $r=\sum_{k=1}^{d} r_k$, $c=\sum_{k=1}^{d} c_k$ and
$$q_{n,d}(L)=q_{n,d}([l_1,\cdots,l_c])=
\sum_{\begin{array}{c}
\scriptstyle J_1\stackrel{\circ}{\cup}\cdots\stackrel{\circ}{\cup}J_{d}=\{1,\dots,c\}\\
\scriptstyle n_k\geq\sum_{i\in J_k}l_i \\
\end{array}}
\prod_{k=1}^{d}{n_k-\sum_{i\in J_k}l_i\choose |J_k|}\cdot|J_k|!\,.$$
\end{theorem}

\begin{remark} \rm
Note that for $d=1$ the above formula coincides with the one due to Robbins:
$$a_{n,1}=\sum_{r=0}^{n-1}\sum_{c=0}^{r}(-1)^{r} 2^{c} 
(n-r-c)!{n-r\choose c}\!\!\!\!\sum_{\begin{array}{l}
\scriptstyle l_1+\cdots+l_{c}=r\\
\scriptstyle l_i\geq 1 
\end{array}}\!\!\!\!\!\!q_{n,1}([l_1,\cdots,l_{c}]),$$
where
\begin{eqnarray*}
\sum_{\begin{array}{l}
\scriptstyle l_1+\cdots+l_{c}=r\\
\scriptstyle l_i\geq 1 
\end{array}}\!\!\!\!\!\!q_{n,1}([l_1,\cdots,l_c])
&=&\!\!\!
\sum_{\begin{array}{l}
\scriptstyle l_1+\cdots+l_{c}=r\\
\scriptstyle l_i\geq 1 
\end{array}}\!\!\!
\sum_{\begin{array}{c}
\scriptstyle J=\{1,\dots,c\}\\
\scriptstyle n\geq r \\
\end{array}}
{n-\sum_{i\in J}l_i\choose |J|}(|J|!)\\
&=&\!\!\!
\sum_{\begin{array}{l}
\scriptstyle l_1+\cdots+l_{c}=r\\
\scriptstyle l_i\geq 1 
\end{array}}{n-r\choose c}c!={r-1\choose c-1}{n-r\choose c}c!\,.
\end{eqnarray*}
\end{remark}

\begin{proof}[Proof of the Theorem]
For $i=1,\dots,n-d$, let $T_{i,d}$ be the set of permutations of $\{1,2,\dots,n\}$
such that $i$ and $i+d$ are {\sl $d$-consecutive}:
$$T_{i,d}=\left\{\pi\in S_n\,:\,|\pi^{-1}(i+d)-\pi^{-1}(i)|=d\right\}\,$$
(for $i=n-d+1,\dots,n$ we consider $T_{i,d}$ as an empty set).
Then by the Inclusion-Exclusion Principle
$$a_{n,d}=\sum_{I\subset \{1,2,\dots,n\}}(-1)^{|I|}\left|\;\bigcap_{i\in I}T_{i,d}\;\right|$$
assuming the convention that when the intersection is made over an empty set of indices then
it is the whole set of permutations $S_n$.
A {\sl $d$-component} of a set of indices $I$ is a maximal subset of $d$-consecutive integers, 
and we denote by $\sharp I$ the number of $d$-components of $I$. 
So the above formula can be rewritten as
$$a_{n,d}=\sum_{r_1=0}^{n_1-1}\sum_{c_1=0}^{r_1}\cdots
\sum_{r_{d}=0}^{n_{d}-1}\sum_{c_{d}=0}^{r_{d}}
(-1)^{r}\;\;\big|\!\!\!\!\!\!\!\!\!\!\!\!\!\!\bigcap_{\begin{array}{c}
\scriptstyle i\in I_1\cup\cdots\cup I_{d}\\
\scriptstyle I_k\subset N_k,\, |I_k|=r_k,\, \sharp I_k=c_k
\end{array}}\!\!\!\!\!\!\!\!\!\!\!\!\!\!\!\!\!\!\!\!T_{i,d}\;\;\;\big|$$
where $N_k=\{1\leq i \leq n: i\equiv k \mod d\}$ and $n_k=|N_k|$.

\begin{remark}\rm
\noindent In order to better illustrate the idea of the proof, which is inspired 
by the one of Robbins in \cite{Robbins}, we give an example of
how a permutation $\pi$ which belongs to the above intersection of $T_{i,d}$'s
can be selected. Assume that $n=12$, $d=2$, $r_1=2$, $c_1=1$, $r_2=3$, $c_2=2$.
$N_1$ and $N_2$ are respectively the odd numbers and the even numbers between $1$ and $12$.
Now we choose $I_1$ and $I_2$: let $l_{i,k}$ be the component sizes 
(each component fixes actually $l_{i,k}+1$ numbers)
and let $j_{i,k}$ be the sizes of the gaps between consecutive components.
Then the choise of $I_1$ and $I_2$ is equivalent to select an integral solution of
$$\left\{\begin{array}{lll}
l_{1,1}=2 &,& l_{i,1}\geq 1 \\
l_{1,2}+l_{2,2}=3 &,& l_{i,2}\geq 1\\
j_{0,1}+j_{1,1}=3 &,& j_{i,1}\geq 0\\
j_{0,2}+j_{1,2}+j_{2,2}=1 &,& j_{i,2}\geq 0
\end{array}\right.\,.$$
For example taking $l_{1,1}=2$, $l_{1,2}=1$, $l_{2,2}=2$, $j_{0,1}=1$, 
$j_{1,1}=2$, $j_{0,2}=0$, $j_{1,2}=1$ and $j_{2,2}=0$, 
we select the set of indices $I_1=\{3,5\}$ and $I_2=\{2,8,10\}$.

\noindent Here is the corresponding table arrangement: 
\begin{center}
\begin{pspicture}(0,-1)(8,3.5)
\def\table{\put(0,0.05){\pspolygon(0,0)(6,0)(6,0.9)(0,0.9)}}
\def\seat{\psellipse(0,0)(0.34,0.34)}
\put(1,0.8){\table
\put(0.5,-.45){\seat \rput(0,0){\bf 2} \rput(0,0.7){$\scriptstyle 2$}}
\put(1.5,-.45){\seat \rput(0,0){4}\rput(0,0.7){$\scriptstyle 4$}}
\put(2.5,-.45){\seat \rput(0,0){6}\rput(0,0.7){$\scriptstyle 6$}}
\put(3.5,-.45){\seat \rput(0,0){\bf 8} \rput(0,0.7){$\scriptstyle 8$}}
\put(4.5,-.45){\seat \rput(0,0){\bf 10} \rput(0,0.7){$\scriptstyle 10$}}
\put(5.5,-.45){\seat \rput(0,0){12}\rput(0,0.7){$\scriptstyle  12$}}
\put(0.1,-.85){\pspolygon[linecolor=gray](0,0)(1.8,0)(1.8,0.8)(0,0.8)}
\put(3.1,-.85){\pspolygon[linecolor=gray](0,0)(2.8,0)(2.8,0.8)(0,0.8)}

\put(0.5,1.45){\seat \rput(0,0){1} \rput(0,-0.7){$\scriptstyle 1$}}
\put(1.5,1.45){\seat \rput(0,0){\bf 3} \rput(0,-0.7){$\scriptstyle 3$}}
\put(2.5,1.45){\seat \rput(0,0){\bf 5} \rput(0,-0.7){$\scriptstyle 5$}}
\put(3.5,1.45){\seat \rput(0,0){7}\rput(0,-0.7){$\scriptstyle 7$}}
\put(4.5,1.45){\seat \rput(0,0){9} \rput(0,-0.7){$\scriptstyle 9$}}
\put(5.5,1.45){\seat \rput(0,0){11} \rput(0,-0.7){$\scriptstyle 11$}}
\put(1.1,1.05){\pspolygon[linecolor=gray](0,0)(2.8,0)(2.8,0.8)(0,0.8)}

\put(0,-1){
\psline{<->}(0,0)(2,0)\uput[d](1,0){$l_{1,2}+1$}
\psline{<->}(2,0)(3,0)\uput[d](2.5,0){$j_{1,2}$}
\psline{<->}(3,0)(6,0)\uput[d](4.5,0){$l_{2,2}+1$}
}

\put(0,2){
\psline{<->}(0,0)(1,0)\uput[u](0.5,0){$j_{0,1}$}
\psline{<->}(1,0)(4,0)\uput[u](3,0){$l_{1,1}+1$}
\psline{<->}(4,0)(6,0)\uput[u](5,0){$j_{1,1}$}
}
}
\end{pspicture}
\end{center}
Now we redistribute the three components selecting a partition
$J_1\stackrel{\circ}{\cup}J_{2}=\{1,2,3\}$, say $J_1=\{1,2\}$ and $J_2=\{3\}$.
This means that the first two components $\{2\}$ and $\{3,5\}$ will go to the odd seats
and the third component $\{8,10\}$ will go to the even seats.
Then we decide the component displacements and orientations: 
for example in the odd seats we place first $\{3,5\}$ reversed and then $\{2\}$
and in the even seats we place $\{8,10\}$ reversed.
To determine the component positions we need the new gap sizes and therefore
 we solve the following two equations
$$\left\{\begin{array}{lll}
j'_{0,1}+j'_{1,1}+j'_{2,1}=1 &,& j'_{i,1}\geq 0 \\
j'_{0,2}+j'_{1,2}=3 &,& j'_{i,2}\geq 0 
\end{array}\right.\,.$$
If we take $j'_{0,1}=0$, $j'_{1,1}=1$, $j'_{2,1}=0$, $j'_{0,2}=1$, $j'_{1,2}=2$ and
we fill the empty places with the remaining numbers $1$, $6$, $9$ and $11$, we obtain the following 
table rearragement:
\begin{center}
\begin{pspicture}(0,-1)(8,3.5)
\def\table{\put(0,0.05){\pspolygon(0,0)(6,0)(6,0.9)(0,0.9)}}
\def\seat{\psellipse(0,0)(0.34,0.34)}
\put(1,0.8){\table
\put(0.5,-.45){\seat \rput(0,0){9} \rput(0,0.7){$\scriptstyle 2$}}
\put(1.5,-.45){\seat \rput(0,0){12}\rput(0,0.7){$\scriptstyle 4$}}
\put(2.5,-.45){\seat \rput(0,0){\bf 10}\rput(0,0.7){$\scriptstyle 6$}}
\put(3.5,-.45){\seat \rput(0,0){\bf 8} \rput(0,0.7){$\scriptstyle 8$}}
\put(4.5,-.45){\seat \rput(0,0){1} \rput(0,0.7){$\scriptstyle 10$}}
\put(5.5,-.45){\seat \rput(0,0){6}\rput(0,0.7){$\scriptstyle  12$}}
\put(1.1,-.85){\pspolygon[linecolor=gray](0,0)(2.8,0)(2.8,0.8)(0,0.8)}

\put(0.5,1.45){\seat \rput(0,0){7} \rput(0,-0.7){$\scriptstyle 1$}}
\put(1.5,1.45){\seat \rput(0,0){\bf 5} \rput(0,-0.7){$\scriptstyle 3$}}
\put(2.5,1.45){\seat \rput(0,0){\bf 3} \rput(0,-0.7){$\scriptstyle 5$}}
\put(3.5,1.45){\seat \rput(0,0){11}\rput(0,-0.7){$\scriptstyle 7$}}
\put(4.5,1.45){\seat \rput(0,0){\bf 2} \rput(0,-0.7){$\scriptstyle 9$}}
\put(5.5,1.45){\seat \rput(0,0){4} \rput(0,-0.7){$\scriptstyle 11$}}
\put(0.1,1.05){\pspolygon[linecolor=gray](0,0)(2.8,0)(2.8,0.8)(0,0.8)}
\put(4.1,1.05){\pspolygon[linecolor=gray](0,0)(1.8,0)(1.8,0.8)(0,0.8)}

\put(0,-1){
\psline{<->}(0,0)(1,0)\uput[d](0.5,0){$j'_{0,2}$}
\psline{<->}(4,0)(6,0)\uput[d](5,0){$j'_{1,2}$}
}

\put(0,2){
\psline{<->}(3,0)(4,0)\uput[u](3.5,0){$j'_{1,1}$}
}
}
\end{pspicture}
\end{center}
\end{remark}

\noindent We go back to the proof. Each $I_k$ is determined by an integral solution of
$$\left\{\begin{array}{lll}
l_{1,k}+l_{2,k}+\cdots+l_{c_k,k}=r_k &,& l_{i,k}\geq 1\\
j_{0,k}+j_{1,k}+\dots +j_{c_k,k}=n_k-r_k-c_k &,& j_{i,k}\geq 0 
\end{array}\right.$$
where $l_{1,k},\dots,l_{c_k,k}$ are the sizes of each component and
$j_{0,k},\dots,j_{c_k,k}$ are the sizes of the gaps.
The counting of the number of integral solutions of this systems yields the
following factor in the formula
$$\prod_{k=1}^{d}{n_k-r_k\choose c_k}\!\!\!\!
\sum_{\begin{array}{c}
\scriptstyle \sum_{i=1}^{c_1}l_{i,1}=r_1\\
\scriptstyle l_{i,1}\geq 1 \\
\end{array}}
\!\!\!\!\!\!\cdots\!\!\!\!\!\!
\sum_{\begin{array}{c}
\scriptstyle \sum_{i=1}^{c_{d}} l_{i,d}=r_{d}\\
\scriptstyle l_{i,d}\geq 1 
\end{array}}.$$
Now, given $I_1,\dots,I_{d}$, we select a permutation 
$\pi\in \left|\bigcap_{i\in I_1\cup\cdots\cup I_{d}}T_{i,d}\right|$:

\begin{itemize}
\item we redistribute the $c=\sum_{k=1}^{d} c_k$ components selecting a partition
$J_1,\dots,J_{d}$ of the set $\{1,\dots,c\}$ (we allow that $J_k$ can be empty),

\item for each set $J_k$ of the partition, we determine the sizes of the new gaps solving
$j'_{0,k}+j'_{1,k}+\dots +j'_{|J_k|,k}=n_k-\sum_{i\in J_k}l_i-|J_k| \quad,\quad j'_{i,k}\geq 0,$

\item we choose the order of components and their orientation in $|J_k|!\cdot 2^{|J_k|}$ ways

\item we fill the empty places with the remaining numbers in $(n-r-c)!$ ways.
\end{itemize}

\noindent Taking into account all these effects we obtain
$$
\sum_{\begin{array}{c}
\scriptstyle J_1\stackrel{\circ}{\cup}\cdots\stackrel{\circ}{\cup}J_{d}=\{1,\dots,c\}\\
\scriptstyle n_k\geq\sum_{i\in J_k}l_i \\
\end{array}}\!\!\!\!\!\!\!(n-r-c)!\prod_{k=1}^{d}{n_k-\sum_{i\in J_k}l_i\choose |J_k|}\cdot|J_k|!\cdot 2^{|J_k|}.$$
Finally, since $c=\sum_{k=1}^{d} |J_k|$, it is easy to get the required formula.
\end{proof}
\end{section}

\begin{section}{Asymptotic analysis: the general case}
\begin{theorem} For $d\geq 0$
$${a_{n,d}\over n!}= e^{-2}\left(1+{4(d-1)\over n}+O\left({1\over n^2}\right)\right).$$
\end{theorem}

\begin{proof} The cases $d=0$ and $d=1$ have been already discussed, so we assume that $d\geq 2$.
By the general formula
$${a_{n,d}\over n!}=\sum {\mbox{polynomial in $n$ of degree $2c$} \over 
\mbox{polynomial in $n$ of degree $r+c$}}.$$
Since $c_k\leq r_k$, the rational functions go to zero faster than $1/n$ as $n$ goes to infinity, 
unless either $r+c=2c$ or $r+c=2c+1$.
In the first case $c_k=r_k$ and $l_{i,k}=1$, in the second case there is one and only one index $k$
such that $c_k+1=r_k$ and $l_{i,k}=2$ for some $i$. By simmetry, we can assume that 
that particular index $k=d$ and multiply the corrisponding term by $d$: 
\begin{eqnarray*}
{a_{n,d}\over n!}&\sim &
\sum_{r_1=0}^{n/d-1}\cdots\sum_{r_{d}=0}^{n/d-1}\left((-2)^{r}{(n-2r)!\over n!}
\prod_{k=1}^{d}{n/d-r_k\choose r_k} 
q_{n,d}([\overbrace{1,\dots,1}^{r}])+\right.\\
&& \!\!\!\!\!\!\!\!\!\!\!\!\!\!\!\!\!\!\!\!\!\!
\left.+d\cdot (-2)^{r-1}{(n-2r+1)!\over n!}
\prod_{k=1}^{d}{n/d-r_k\choose r_k} {n/d-r_{d}\choose r_{d}-1} (r_{d}-1)
q_{n,d}([\overbrace{1,\dots,1}^{r-2},2])\right),
\end{eqnarray*}
that is
\begin{eqnarray*}
{a_{n,d}\over n!}&\sim &
\sum_{r_1=0}^{n/d-1}\cdots\sum_{r_{d}=0}^{n/d-1}{q_{n,d}([\overbrace{1,\dots,1}^{r}])\over n^{\underline{2r}}}
\prod_{k=1}^{d}{(-2)^{r_k}\over r_k!}\prod_{k=1}^{d}{(n/d)^{\underline{2r_k}}\over  (n/d)^{\underline{r_k}}}+\\
&& \!\!\!\!\!\!\!\!\!\!\!\!\!\!\!\!\!\!\!\!\!\!
+2d\cdot \sum_{r_1=0}^{n/d-1}\cdots\sum_{r_{d-1}=0}^{n/d-1}\sum_{r_{d}=2}^{n/d-1}
{q_{n,d}([\overbrace{1,\dots,1}^{r-2},2])\over n^{\underline{2r-1}}}
\prod_{k=1}^{d-1}{(-2)^{r_k}\over r_k!}\cdot{(-2)^{r_{d}-2}\over (r_{d}-2)!}
\prod_{k=1}^{d-1}{(n/d)^{\underline{2r_k}}\over (n/d)^{\underline{r_k}}}\cdot
{(n/d)^{\underline{2r_{d}-1}}\over  (n/d)^{\underline{r_{d}}}}.
\end{eqnarray*}
We start considering the second term. Since
\begin{eqnarray*}
q_{n,d}([\overbrace{1,\dots,1}^{r-2},2])&\sim&\!\!\!\!\!\!
\sum_{\begin{array}{c}
\scriptstyle \sum_{k=1}^{d}r'_{k}=r-1\\
\scriptstyle r'_{k}\geq 0 \\
\end{array}} {r-1 \choose r'_1,\dots, r'_{d}}
\prod_{k=1}^{d-1}{n/d-r'_k\choose r'_k}\cdot r'_k!\cdot {n/d-r'_d-1\choose r'_d}\cdot r'_d!\\
&\sim&\!\!\!\!\!\!
\sum_{\begin{array}{c}
\scriptstyle \sum_{k=1}^{d}r'_{k}=r-1\\
\scriptstyle r'_{k}\geq 0 \\
\end{array}} {r-1 \choose r'_1,\dots, r'_{d}}
\prod_{k=1}^{d-1}{(n/d)^{\underline{2r'_k}}\over (n/d)^{\underline{r'_k}}} \cdot
{(n/d)^{\underline{2r'_d+1}}\over (n/d)^{\underline{r'_d+1}}}\\
&\sim&\left({n\over d}\right)^{r-1}\!\!\!\!\!\!
\sum_{\begin{array}{c}
\scriptstyle \sum_{k=1}^{d}r'_{k}=r-1\\
\scriptstyle r'_{k}\geq 0 \\
\end{array}} {r-1 \choose r'_1,\dots, r'_{d}}=
\left({n\over d}\right)^{r-1}d^{r-1}=n^{r-1}
\end{eqnarray*}
then it becomes
$$
2\cdot \sum_{r_1=0}^{n/d-1}\cdots\sum_{r_{d-1}=0}^{n/d-1}\sum_{r_{d}=2}^{n/d-1}
{n^{r-1}\over n^{2r-1}}\prod_{k=1}^{d-1}{(-2/d)^{r_k}\over r_k!}\cdot{(-2/d)^{r_{d}-2}\over (r_{d}-2)!}\cdot n^{r-1}\sim
2\left(e^{-2/d}\right)^d{1 \over n}={2 e^{-2}\over n}.
$$
Now we consider the first term. Since
\begin{eqnarray*}
q_{n,d}([\overbrace{1,\dots,1}^{r}])&\sim&\!\!\!\!\!\!
\sum_{\begin{array}{c}
\scriptstyle \sum_{k=1}^{d}r'_{k}=r\\
\scriptstyle r'_{k}\geq 0 \\
\end{array}} {r \choose r'_1,\dots, r'_{d}}
\prod_{k=1}^{d}{n/d-r'_k\choose r'_k}\cdot r'_k!\\
&\sim&\!\!\!\!\!\!
\sum_{\begin{array}{c}
\scriptstyle \sum_{k=1}^{d}r'_{k}=r\\
\scriptstyle r'_{k}\geq 0 \\
\end{array}} {r \choose r'_1,\dots, r'_{d}}
\prod_{k=1}^{d}{(n/d)^{\underline{2r_k}}\over (n/d)^{\underline{r_k}}}\\
&\sim&
\left({n\over d}\right)^r\!\!\!\!\!\!
\sum_{\begin{array}{c}
\scriptstyle \sum_{k=1}^{d}r'_{k}=r\\
\scriptstyle r'_{k}\geq 0 \\
\end{array}} {r \choose r'_1,\dots, r'_{d}}
\prod_{k=1}^{d}{1-{1\over 2}(2r'_k)(2r'_k-1)\cdot{d\over n}\over 1-{1\over 2}r'_k(r'_k-1)\cdot{d\over n}}\\
&\sim&
\left({n\over d}\right)^r\!\!\!\!\!\!
\sum_{\begin{array}{c}
\scriptstyle \sum_{k=1}^{d}r'_{k}=r\\
\scriptstyle r'_{k}\geq 0 \\
\end{array}} {r \choose r'_1,\dots, r'_{d}}
\left(1-\left({3\over 2}\sum_{k=1}^{d}r'_k(r'_k-1)+\sum_{k=1}^{d}r'_k\right)\cdot{d\over n}\right)\\
&\sim&
\left({n\over d}\right)^r
\left(d^r-\left({3\over 2}dr(r-1)d^{r-2}+drd^{r-1}\right)\cdot{d\over n}\right)\\
&\sim& n^r
\left(1-\left({3\over 2}r(r-1)+dr\right){1\over n}\right)
\end{eqnarray*}
then it becomes
$$\sum_{r_1=0}^{n/d-1}\cdots\sum_{r_{d}=0}^{n/d-1}
{1-\left({3\over 2}r(r-1)+dr\right)\cdot{1\over n}\over 1-r(2r-1)\cdot{1\over n}}
\prod_{k=1}^{d}{(-2/d)^{r_k}\over r_k!}\prod_{k=1}^{d}
\left({1-r_k(2r_k-1)\cdot{d\over n}\over 1-{1\over 2}r_k(r_k-1)\cdot{d\over n}}\right)
$$
that is
$$\sum_{r_1=0}^{n/d-1}\cdots\sum_{r_{d}=0}^{n/d-1}
\left(1+\left({r^2+r\over 2}-dr\right)\cdot{1\over n}\right)
\prod_{k=1}^{d}{(-2/d)^{r_k}\over r_k!}\prod_{k=1}^{d}
\left(1-\left({3r_k(r_k-1)\over 2}+r_k\right)\cdot {d\over n}\right).
$$
Riminding that $r=\sum_{k=1}^d r_k$, it is equivalent to
\begin{eqnarray*}e^{-2}+{1\over n}\cdot\left.
\sum_{r_1=0}^{n/d-1}\cdots\sum_{r_{d}=0}^{n/d-1}
\left({1-3d\over 2}\sum_{k=1}^{d} r_k(r_{k}-1)\prod_{k=1}^{d}{(-2/d)^{r_k}\over r_k!}\right.\right.+\\
+(1-2d)\sum_{k=1}^{d} r_k\prod_{k=1}^{d}{(-2/d)^{r_k}\over r_k!}+
\left.\left.\sum_{k=1}^{d}\sum_{k'=k+1}^{d} r_k r_{k'}\prod_{k=1}^{d}{(-2/d)^{r_k}\over r_k!}\right)\right.
\end{eqnarray*}
and taking the limit of the remaining sums we obtain
$$e^{-2}\left(1+{1\over n}\cdot\left({1-3d\over 2}\left({-2\over d}\right)^2 d +
(1-2d)\left({-2\over d}\right)d+\left({-2\over d}\right)^2{d(d-1)\over 2}\right)\right)$$
that is
$$e^{-2}\left(1+{4d-6\over n}\right).$$
Finally, putting all together we find that
$${a_{n,d}\over n!}\sim e^{-2}\left(1+{4d-6\over n}\right)+{2 e^{-2}\over n}=e^{-2}\left(1+{4(d-1)\over n}\right).$$
\end{proof}
\end{section}

\end{document}